\documentclass[english]{iopart}
\usepackage[T1]{fontenc}
\usepackage[latin9]{inputenc}
\usepackage{color}
\usepackage{float}
\usepackage{booktabs}
\usepackage{amsbsy}
\usepackage{amstext}
\usepackage{graphicx}
\usepackage{amssymb}

\makeatletter

\providecommand{\tabularnewline}{\\}

\usepackage{iopams}
\usepackage{setstack}

\@ifundefined{showcaptionsetup}{}{%
 \PassOptionsToPackage{caption=false}{subfig}}
\usepackage{subfig}
\makeatother

\usepackage{babel}

\begin{document}
\title[Perturbation Resilience and Superiorization]{Perturbation Resilience and Superiorization of Iterative Algorithms}

\author{Y Censor$^{1}$, R Davidi$^{2}$ and G T Herman$^{2}$}

\address{$^{1}$Department of Mathematics, University of Haifa, Mount Carmel,
Haifa 31905, Israel\\
$^{2}$Department of Computer Science, Graduate Center, City University
of New York, New York, NY 10016, USA}

\ead{$^{1}$yair@math.haifa.ac.il}
\begin{abstract}
\noindent Iterative algorithms aimed at solving some problems are
discussed. For certain problems, such as finding a common point in
the intersection of a finite number of convex sets, there often exist
iterative algorithms that impose very little demand on computer resources.
For other problems, such as finding that point in the intersection
at which the value of a given function is optimal, algorithms tend
to need more computer memory and longer execution time. A methodology
is presented whose aim is to produce automatically for an iterative
algorithm of the first kind a {}``superiorized version'' of it that
retains its computational efficiency but nevertheless goes a long
way towards solving an optimization problem. This is possible to do
if the original algorithm is {}``perturbation resilient,'' which
is shown to be the case for various projection algorithms for solving
the consistent convex feasibility problem. The superiorized versions
of such algorithms use perturbations that drive the process in the
direction of the optimizer of the given function. After presenting
these intuitive ideas in a precise mathematical form, they are illustrated
in image reconstruction from projections for two different projection
algorithms superiorized for the function whose value is the total
variation of the image.
\end{abstract}



\submitto{Inverse Problems}

\maketitle
\maketitle

\section{Introduction}

\textit{\emph{We first motivate and describe }}our ideas in a not
fully general context, in which \textit{\emph{superiorization}}\emph{
}is envisioned as lying in-between the methodologies of optimization
and of feasibility seeking. With a feasible solution one settles for
a point that just fulfills a set of constraints, whereas solving a
constrained optimization problem calls for finding a feasible point
that optimizes a given objective function. Generally speaking, optimization
is logically and computationally a more demanding task than that of
finding just any feasible point. We show that, without employing an
optimization algorithm, it is possible to use certain iterative methods,
designed for (the less demanding) feasibility problems, in a way that
will steer the iterates toward a point that is \textit{superior,}
but not necessarily optimal, in a well-defined sense. The advantage
of superiorization is that it allows us to solve significant problems
by using powerful feasibility seeking methods, see, e.g., \cite{cccdh09}
and references therein, and reach a \textit{superior feasible point}
without resorting to optimization techniques. We now explain this
with more details.

Many significant real-world problems are modeled by constraints that
force the sought-after solution point to fulfill conditions imposed
by the physical nature of the problem. Such a modeling approach often
leads to a \textit{convex feasibility problem }of the form\begin{equation}
\text{find }\boldsymbol{x}^{\ast}\in C=\bigcap_{i=1}^{I}C_{i},\label{eq:cfp}\end{equation}
where the sets $C_{i}\subseteq\mathbb{R}^{J}$ are closed convex subsets
of the Euclidean space $\mathbb{R}^{J}$, see \cite{bb96,byrnebook,chinneck-book}
or \cite[Chapter 5]{CZ97} for this broad topic. In many real-world
problems the underlying system is very large (huge values of $I$
and $J$) and often very sparse. In these circumstances \textit{projection
methods} have proved to be effective. They are iterative algorithms
that use projections onto sets while relying on the general principle
that when a family of closed and convex sets is present, then projections
onto the individual sets are easier to perform than projections onto
other sets, such as their intersection as in (\ref{eq:cfp}), that
are derived from them.

Projection methods can have various algorithmic structures (some of
which are particularly suitable for parallel computing) and they also
possess desirable convergence properties and good initial behavior
patterns \cite{bb96,CZ97,c96,Jamo97,Imag97,H09,Pier84}. The main
advantage of projection methods, which makes them successful in real-world
applications, is computational. They commonly have the ability to
handle huge-size problems of dimensions beyond which more sophisticated
methods cease to be efficient or even applicable due to memory requirements.
(For a justification of this claim see the various examples provided
in \cite{cccdh09}.) This is so because the building bricks of a projection
algorithm (which are the projections onto the given individual sets)
are easy to perform, and because the algorithmic structure is either
sequential or simultaneous, or in-between, as in the block-iterative
projection methods or in the more recently invented string-averaging
projection methods. The number of sets used simultaneously in each
iteration in block-iterative methods and the number and lengths of
strings used in each iteration in string-averaging methods are variable,
which provides great flexibility in matching the implementation of
the algorithm with the parallel architecture at hand; for block-iterative
methods see, e.g., \textcolor{black}{\cite{ac89,Baus06,Butn90,cgg01,Imag97,dhc09,eGG81,Gonz01,Lopu97,Otta88,pen09}
and for string-averaging methods see, e.g., \cite{bmr03,bdhk07,ceh01,cs09,ct03,crombez02,pen09,rh03}.}

The key to superiorization is our recent discovery \cite{bdhk07,dhc09,hd08}
that two principal prototypical algorithmic schemes of projection
methods: string-averaging projections (SAP) and block-iterative projections
(BIP), which include as special cases a variety of projection methods
for the convex feasibility problem, are\emph{ }\textit{\emph{bounded
perturbations resilient}} in the sense that the convergence of sequences
generated by them continues to hold even if the iterates are perturbed
in every iteration. We harness this resilience to bounded perturbations
to steer the iterates to not just any feasible point but to a superior
(in a well-defined sense) feasible point of (\ref{eq:cfp}).

Our motivation is the desire to create a new methodology that will
significantly improve methods for the solution of inverse problems
in image reconstruction from projections, intensity-modulated radiation/proton
therapy (IMRT/IMPT) and in other real-world problems such as electron
microscopy (EM). Our work \cite{bdhk07,dhc09}, as well as the examples
given below, indicate that our objective is achievable and show how
algorithms can incorporate\emph{ }\textit{\emph{perturbations}} in
order to perform \textit{\emph{superiorization}}. 

The superiorization methodology has in fact broader applicability
than what has been discussed until now and its mathematical specification
in the next section reflects this. However, all our specific examples
will be chosen from the field that we used as our motivation in this
introductory section.

\section{Specification of the superiorization methodology\label{sec:methodolgy}}

The superiorization principle relies on the bounded perturbation resilience
of algorithms. Therefore we define this notion next in a general setting
within $\mathbb{R}^{J}$. 

We introduce the notion of a \emph{problem} \emph{structure} $\left\langle \mathbb{T},\mathcal{P}r\right\rangle $,
where $\mathbb{T}$ is a nonempty \emph{problem set} and $\mathcal{P}r$
is a function on $\mathbb{T}$ such that, for all $T\in\mathbb{T}$,
$\mathcal{P}r_{T}:\mathbb{R}^{J}\rightarrow\mathbb{R}_{+}$, where
$\mathbb{R}_{+}$ is the set of nonnegative real numbers. Intuitively
we think of $\mathcal{P}r_{T}\left(\boldsymbol{x}\right)$ as a measure
of how {}``far'' $\boldsymbol{x}$ is from being a solution of $T$.
In fact, we call $\boldsymbol{x}$ a \emph{solution} of $T$ if $\mathcal{P}r_{T}\left(\boldsymbol{x}\right)=0$.

For example, for the convex feasibility problem (\ref{eq:cfp})\begin{equation}
\begin{array}{cl}
\mathbb{T}=\left\{ \left\{ C_{1},\ldots,C_{I}\right\} \right| & I\mbox{\,\ is a positive integer and, for \ensuremath{1\leq i\leq I,}}\\
 & \left.C_{i}\mbox{\ is a closed convex subset of \ensuremath{\mathbb{R}^{J}}}\right\} \end{array}\label{eq:EX-CFP}\end{equation}
 and \begin{equation}
\mathcal{P}r_{\left\{ C_{1},\ldots,C_{I}\right\} }\left(\boldsymbol{x}\right)=\sqrt{\sum_{i=1}^{I}\left(d\left(\boldsymbol{x},C_{i}\right)\right)^{2}},\label{eq:Res}\end{equation}
where $d\left(\boldsymbol{x},C_{i}\right)$ is the Euclidean distance
of $\boldsymbol{x}$ from the set $C_{i}$. Clearly, in this case
$\boldsymbol{x}$ is a solution of $\left\{ C_{1},\ldots,C_{I}\right\} $
as defined in the previous paragraph if, and only if, $\boldsymbol{x}\in C$
as defined in (\ref{eq:cfp}).\textbf{ }\\
\textbf{}\\
\textbf{Definition 1.\label{def:resilient} }An \emph{algorithm
}\textbf{$\mathbf{P}$ }for $\left\langle \mathbb{T},\mathcal{P}r\right\rangle $
assigns to each $T\in\mathbb{T}$ an algorithmic operator $\mathbf{P}_{T}:\mathbb{R}^{J}\rightarrow\mathbb{R}^{J}$.
$\mathbf{P}$ is said to be \emph{bounded perturbations resilient
}if, for all $T\in\mathbb{T}$, the following is the case: if the
sequence $\left\{ \left(\mathbf{P}_{T}\right)^{k}\boldsymbol{x}\right\} _{k=0}^{\infty}$
converges to a solution of $T$ for all $\boldsymbol{x}\in\mathbb{R}^{J}$,
then any sequence $\left\{ \boldsymbol{x}^{k}\right\} _{k=0}^{\infty}$
of points in $\mathbb{R}^{J}$ also converges to a solution of $T$
provided that, for all $k\geq0$, \begin{equation}
\boldsymbol{x}^{k+1}=\mathbf{P}_{T}\left(\boldsymbol{x}^{k}+\beta_{k}\boldsymbol{v}^{k}\right),\label{eq:resilient}\end{equation}
where $\beta_{k}\boldsymbol{v}^{k}$ are bounded perturbations, meaning
that $\beta_{k}$ are real nonnegative numbers such that ${\displaystyle \sum\limits _{k=0}^{\infty}}\beta_{k}\,<\infty$
and the sequence $\left\{ \boldsymbol{v}^{k}\right\} _{k=0}^{\infty}$
is bounded. 

We give next specific instances of bounded perturbations resilient
algorithms for solving the convex feasibility problem as in (\ref{eq:EX-CFP})
and (\ref{eq:Res}), from the classes of SAP and BIP methods. We do
this by defining $\mathbf{P}_{\left\{ C_{1},\ldots,C_{I}\right\} }$
for an arbitrary but fixed element $\left\{ C_{1},\ldots,C_{I}\right\} $
of $\mathbb{T}$ of (\ref{eq:EX-CFP}) for the different algorithms
$\mathbf{P}$. For any nonempty closed convex subset $M$ of $\mathbb{R}^{J}$
and any $\boldsymbol{x}\in\mathbb{R}^{J},$ the orthogonal projection
of $\boldsymbol{x}$ onto $M$ is the point in $M$ that is nearest
(by the Euclidean distance) to $\boldsymbol{x}$; it is denoted by
$P_{M}\boldsymbol{x}$.

To define $\mathbf{P}_{\left\{ C_{1},\ldots,C_{I}\right\} }$ for
the SAP instances, we make use of \textit{index vector}\textit{\emph{s}}\textit{,
}\textit{\emph{which are}} nonempty ordered sets $t=\left(t_{1},\dots,t_{N}\right)$,
where $N$ is an arbitrary positive integer, whose elements $t_{n}$
are in the set $\left\{ 1,...,I\right\} .$ For an index vector $t$
we define the composite operator \begin{equation}
P\left[t\right]=P_{C_{t_{N}}}\cdots P_{C_{t_{1}}}.\label{eq:P_C[t]}\end{equation}
A finite set $\Omega$ of index vectors is called \textit{fit} if,
for each $i\in\left\{ 1,...,I\right\} $, there exists $t=\left(t_{1},\dots,t_{N}\right)\in\Omega$
such that $t_{n}=i$ for some $n\in\left\{ 1,...,N\right\} .$ If
$\Omega$ is a fit set of index vectors, then a function $\omega:\Omega\rightarrow\mathbb{R}_{++}=\left(0,\infty\right)$
is called a \textit{fit weight function} if $\sum_{t\in\Omega}\omega\left(t\right)=1.$
A pair $\left(\Omega,\omega\right)$ consisting of a fit set of index
vectors and a fit weight function defined on it was called an \textit{amalgamator}
in \cite{bdhk07}. For each amalgamator $\left(\Omega,\omega\right),$
we define the algorithmic operator $\mathbf{P}_{\left\{ C_{1},\ldots,C_{I}\right\} }:\mathbb{R}^{J}\rightarrow\mathbb{R}^{J}$
by\begin{equation}
\mathbf{P}_{\left\{ C_{1},\ldots,C_{I}\right\} }\boldsymbol{x}\mathbf{=}\sum_{t\in\Omega}\omega\left(t\right)P\left[t\right]\boldsymbol{x}.\label{eq: Bold p}\end{equation}
For this algorithmic operator we have the following \textit{\emph{bounded
perturbations resilience theorem.}}\textbf{\textit{\emph{ }}}\\
\textbf{\textit{\emph{}}}\\
\textbf{\textit{\emph{Theorem 1.}}}\label{Theorem-1}\textit{ }\textbf{\cite[Section II]{bdhk07}}\textit{
If} $C$ \emph{of (\ref{eq:cfp}) is nonempty, }$\left\{ \beta_{k}\right\} _{k=0}^{\infty}$\textit{
is a sequence of nonnegative real numbers such that} $\sum_{k=0}^{\infty}\beta_{k}<\infty$
\textit{and }$\left\{ \boldsymbol{v}^{k}\right\} _{k=0}^{\infty}$
\textit{is a bounded sequence of points in} $\mathbb{R}^{J},$ \textit{then
for any amalgamator} $\left(\Omega,\omega\right)$ \textit{and any
$\boldsymbol{x}^{0}\in\mathbb{R}^{J}$, the sequence }$\left\{ \boldsymbol{x}^{k}\right\} _{k=0}^{\infty}$\textit{
generated by\begin{equation}
\boldsymbol{x}^{k+1}=\mathbf{P}_{\left\{ C_{1},\ldots,C_{I}\right\} }\left(\boldsymbol{x}^{k}+\beta_{k}\boldsymbol{v}^{k}\right),\text{ }\forall k\geq0,\label{eq:Theorem amalgamated}\end{equation}
converges, and its limit is in} $C.$ (The statement of this theorem
in \cite{bdhk07} is for positive $\beta_{k}$s, but the proof given
there applies to nonnegative $\beta_{k}$s.)\\
\\
\textbf{Corollary 1.} \emph{For any amalgamator }$(\Omega,\omega)$,
\emph{the algorithm $\mathbf{P}$ defined by the algorithmic operator
$\mathbf{P}_{\left\{ C_{1},\ldots,C_{I}\right\} }$ is bounded perturbations
resilient.}\\
\textbf{Proof. }Assume that for $T=\left\{ C_{1},\ldots,C_{I}\right\} $
the sequence $\left\{ \left(\mathbf{P}_{T}\right)^{k}\boldsymbol{x}\right\} _{k=0}^{\infty}$
converges to a solution of $T$ for all $\boldsymbol{x}\in\mathbb{R}^{J}$.
This implies, in particular, that $C$ of (\ref{eq:cfp}) is nonempty.
By Definition 1, we need to show that any sequence $\left\{ \boldsymbol{x}^{k}\right\} _{k=0}^{\infty}$
of points in $\mathbb{R}^{J}$ also converges to a solution of $T$
provided that, for all $k\geq0$, (\ref{eq:resilient}) is satisfied
when the $\beta_{k}\boldsymbol{v}^{k}$ are \textit{\emph{bounded}}\textit{
}\textit{\emph{perturbations}}\emph{.} Under our assumptions, this
follows from Theorem 1. $\square$\\

\textcolor{black}{Next we look at a member of the family of BIP}\textcolor{black}{\emph{
}}\textcolor{black}{methods}\emph{.} Considering the convex feasibility
problem (\ref{eq:cfp}), for $1\leq u\leq U,$ let $B_{u}$ be an
ordered set $\left(b_{u,1},\dots,b_{u,\left|B_{u}\right|}\right)$
of elements of $\left\{ 1,\dots,I\right\} $ ($\left|B_{u}\right|$
denotes the cardinality of $B_{u}$). We \textcolor{black}{call such
a $B_{u}$ a }\textcolor{black}{\emph{block}}\textcolor{black}{{} $\;$and
define the (composite) algorithmic operator $\mathbf{\mathbf{Q}}_{\left\{ C_{1},\ldots,C_{I}\right\} }:\mathbb{R}^{J}\rightarrow\mathbb{R}^{J}$
by \begin{equation}
\mathbf{\mathbf{Q}}_{\left\{ C_{1},\ldots,C_{I}\right\} }=Q_{U}\cdots Q_{1},\label{eq: Blod Q operator}\end{equation}
where, for $\boldsymbol{x}\in\mathbb{R}^{J}$ and $1\leq u\leq U$,
\begin{equation}
Q_{u}\boldsymbol{x}=\frac{1}{R}\sum_{i\in B_{u}}P_{C_{i}}\boldsymbol{x}+\frac{R-\left|B_{u}\right|}{R}\boldsymbol{x},\label{eq: Q_u operator}\end{equation}
and \begin{equation}
R=\max\left\{ \left|B_{u}\right|\mid1\leq u\leq U\right\} .\label{eq:definition of R}\end{equation}
The iterative procedure $\boldsymbol{x}^{k+1}=\mathbf{\mathbf{Q}}_{\left\{ C_{1},\ldots,C_{I}\right\} }\boldsymbol{x}^{k}$
is a member of the family of BIP}\textcolor{black}{\emph{ }}\textcolor{black}{methods.
}For this algorithmic operator we have the following \textit{\emph{bounded
perturbations resilience theorem.}}\textbf{\textit{\emph{ }}}\\
\textbf{\textit{\emph{}}}\\
\textbf{\textit{\emph{Theorem 2.}}}\label{Theorem 2}\textit{ }\textbf{\cite{dhc09}}\textit{
If}\textit{\emph{ $C$}}\emph{ of (\ref{eq:cfp}) is}\textit{\emph{
}}\textit{nonempty,}\textit{\emph{ $\left\{ 1,\dots,I\right\} =\bigcup_{u=1}^{U}B_{u}$,}}\emph{
$\left\{ \beta_{k}\right\} _{k=0}^{\infty}$}\textit{\emph{ }}\textit{is
a sequence of nonnegative real numbers such that} $\sum_{k=0}^{\infty}\beta_{k}<\infty$,
$\left\{ \boldsymbol{v}^{k}\right\} _{k=0}^{\infty}$ \textit{be a
bounded sequence of points in} $\mathbb{R}^{J},$ \textit{then for
any $\boldsymbol{x}^{0}\in\mathbb{R}^{J}$, the sequence }$\left\{ \boldsymbol{x}^{k}\right\} _{k=0}^{\infty}$\textit{
generated by\begin{equation}
\boldsymbol{x}^{k+1}=\mathbf{\mathbf{Q}}_{\left\{ C_{1},\ldots,C_{I}\right\} }\left(\boldsymbol{x}^{k}+\beta_{k}\boldsymbol{v}^{k}\right),\text{ }\forall k\geq0,\label{eq:Theorem Block procedure}\end{equation}
converges, and its limit is in} $C.$ (This is a special case of Theorem
2 in \cite{dhc09} given here without a relaxation parameter. Also,
that theorem is stated for positive $\beta_{k}$s, but the proof given
there applies to nonnegative $\beta_{k}$s.)\\
\\
\textbf{Corollary 2.} \emph{The algorithm $\mathbf{Q}$ defined
by the algorithmic operator $\mathbf{Q}_{\left\{ C_{1},\ldots,C_{I}\right\} }$
is bounded perturbations resilient.}\\
\textbf{Proof.} Replace in the proof of Corollary 1 \emph{$\mathbf{P}$
}by\emph{ $\mathbf{Q}$} and Theorem 1 by Theorem 2. $\square$\\

Further \textit{\emph{bounded perturbations resilience theorems are
available in a Banach space setting, see \cite{brz05,brz08}. Thus
the theory of bounded perturbations resilient algorithms already contains
some solid mathematical results. As opposed to this, the}} superiorization
theory that we present next is at the stage of being a collection
of heuristic ideas, a full mathematical theory still needs to be developed.
However, there are practical demonstrations of its potential usefulness;
see \cite{bdhk07,dhc09,hd08} and the illustrations in Section \ref{sec:tv}
below.

For a problem structure $\left\langle \mathbb{T},\mathcal{P}r\right\rangle $,
$T\in\mathbb{T}$, $\varepsilon\in\mathbb{R}_{++}$ and a sequence
$S=\left\{ \boldsymbol{x}^{k}\right\} _{k=0}^{\infty}$ of points
in $\mathbb{R}^{J}$, we use $O\left(T,\varepsilon,S\right)$ to denote
the $\boldsymbol{x}\in\mathbb{R}^{J}$ that has the the following
properties: $\mathcal{P}r_{T}(\boldsymbol{x})\leq\varepsilon$\emph{
}and there is a nonnegative integer $K$ such that $\boldsymbol{x}^{K}=\boldsymbol{x}$
and, for all nonnegative integers $\ell<K$, $\mathcal{P}r_{T}\left(\boldsymbol{x}^{\ell}\right)>\varepsilon$.
Clearly, if there is such an $\boldsymbol{x}$, then it is unique.
If there is no such $\boldsymbol{x}$, then we say that $O\left(T,\varepsilon,S\right)$
is undefined. The intuition behind this definition is the following:
if we think of $S$ as the (infinite) sequence of points that is produced
by an algorithm (intended for the problem $T$) without a termination
criterion, then $O\left(T,\varepsilon,S\right)$ is the output produced
by that algorithm when we add to it instructions that makes it terminate
as soon as it reaches a point at which the value of $\mathcal{P}r_{T}$
is not greater than $\varepsilon$. The following result is obvious.\textbf{\textit{\emph{}}}\\
\textbf{\textit{\emph{}}}\\
\textbf{\textit{\emph{Lemma 1.}}}\label{lem:If--is}\emph{ If $\mathcal{P}r_{T}$
is continuous and the sequence $S$ converges to a solution of $T$,
then $O\left(T,\varepsilon,S\right)$ is defined and $\mathcal{P}r_{T}\left(O\left(T,\varepsilon,S\right)\right)\leq\varepsilon$.}\\

Given an algorithm\emph{ }\textbf{$\mathbf{P}$ }for a problem structure
$\left\langle \mathbb{T},\mathcal{P}r\right\rangle $, a $T\in\mathbb{T}$
and an $\bar{\boldsymbol{x}}\in\mathbb{R}^{J}$, let $R\left(T,\bar{\boldsymbol{x}}\right)=\left\{ \left(\mathbf{P}_{T}\right)^{k}\boldsymbol{x}\right\} _{k=0}^{\infty}$.
For a function $\phi:\mathbb{R}^{J}\rightarrow\mathbb{R}$, the \emph{superiorization
methodology }should provide us with an algorithm that produces a sequence
$S\left(T,\bar{\boldsymbol{x}},\phi\right)=\left\{ \boldsymbol{x}^{k}\right\} _{k=0}^{\infty}$,
such that for any $\varepsilon\in\mathbb{R}_{++}$ and $\bar{\boldsymbol{x}}\in\mathbb{R}^{J}$
for which $\mathcal{P}r_{T}\left(\bar{\boldsymbol{x}}\right)>\varepsilon$
and $O\left(T,\varepsilon,R\left(T,\bar{\boldsymbol{x}}\right)\right)$
is defined, $O\left(T,\varepsilon,S\left(T,\bar{\boldsymbol{x}},\phi\right)\right)$
is also defined and $\phi\left(O\left(T,\varepsilon,S\left(T,\bar{\boldsymbol{x}},\phi\right)\right)\right)<\phi\left(O\left(T,\varepsilon,R\left(T,\bar{\boldsymbol{x}}\right)\right)\right)$.
This is of course too ambitious in its full generality and so here
we analyze only a special case, but one that is still quite general.
We now list our assumptions for the special case for which we discuss
details of the superiorization methodology.\\
\\
\textbf{\textit{\emph{Assumptions}}}
\begin{enumerate}
\item $\left\langle \mathbb{T},\mathcal{P}r\right\rangle $ is a problem
structure such that $\mathcal{P}r_{T}$ is continuous for all $T\in\mathbb{T}$.
\item $\mathbf{P}$ is a bounded perturbation resilient algorithm for $\left\langle \mathbb{T},\mathcal{P}r\right\rangle $
such that, for all $T\in\mathbb{T}$, $\mathbf{P}_{T}$ is continuous
and, if $\boldsymbol{x}$ is not a solution of $T$, then $\mathcal{P}r_{T}\left(\mathbf{P}_{T}\boldsymbol{x}\right))<\mathcal{P}r_{T}\left(\boldsymbol{x}\right)$.
\item $\phi$ is a convex function.\\

\end{enumerate}
We now describe, under these assumptions, the algorithm to produce
the sequence $S\left(T,\bar{\boldsymbol{x}},\phi\right)=\left\{ \boldsymbol{x}^{k}\right\} _{k=0}^{\infty}$.

The algorithm assumes that we have available a summable sequence $\left\{ \gamma_{\ell}\right\} _{\ell=0}^{\infty}$
of positive real numbers. It is easy to generate such sequences; e.g.,
we can use $\gamma_{\ell}=a^{\ell}$, where $0<a<1$. The algorithm
generates, simultaneously with the sequence $\left\{ \boldsymbol{x}^{k}\right\} _{k=0}^{\infty}$,
sequences $\left\{ \boldsymbol{v}^{k}\right\} _{k=0}^{\infty}$ and
$\left\{ \beta_{k}\right\} _{k=0}^{\infty}$. The latter will be generated
as a subsequence of $\left\{ \gamma_{\ell}\right\} _{\ell=0}^{\infty}$.
Clearly, the resulting sequence $\{\beta_{k}\}_{k=0}^{\infty}$ of
positive real numbers will be summable. We first specify the algorithm
and then discuss it. The algorithm depends on the specified $\bar{\boldsymbol{x}}$,
$\phi$, $\left\{ \gamma_{\ell}\right\} _{\ell=0}^{\infty}$, $\mathcal{P}r_{T}$
and $\mathbf{P}_{T}$. It makes use of a logical variable called \emph{continue
}and also of the concept of a subgradient of the convex function $\phi$.\\
\\
\textbf{\textit{\emph{Superiorized Version of Algorithm $\mathbf{P}$}}}
\begin{enumerate}
\item \textbf{set} $k=0$
\item \textbf{set} $\boldsymbol{x}^{k}=\bar{\boldsymbol{x}}$
\item \textbf{set} $\ell=0$
\item \textbf{repeat}
\item $\qquad$\textbf{set $\boldsymbol{g}$ }to a subgradient of $\phi$
at $\boldsymbol{x}^{k}$
\item $\qquad$\textbf{if $\left\Vert \boldsymbol{g}\right\Vert >0$}
\item $\qquad\qquad$\textbf{then set $\boldsymbol{v}^{k}=-\boldsymbol{g}/\left\Vert \boldsymbol{g}\right\Vert $}
\item $\qquad\qquad$\textbf{else set $\boldsymbol{v}^{k}=\boldsymbol{g}$}
\item \textbf{$\qquad$set} \emph{continue = true}
\item $\qquad$\textbf{while} \emph{continue}
\item $\qquad\qquad$\textbf{set} $\beta_{k}=\gamma_{\ell}$
\item $\qquad\qquad$\textbf{set} $\boldsymbol{y}=\boldsymbol{x}^{k}+\beta_{k}\boldsymbol{v}^{k}$
\item $\qquad\qquad$\textbf{if $\phi\left(\boldsymbol{y}\right)\leq\phi\left(\boldsymbol{x}^{k}\right)$
and $\mathcal{P}r_{T}\left(\mbox{\ensuremath{\mathbf{P}_{T}\boldsymbol{y}}}\right)<\mathcal{P}r_{T}\left(\boldsymbol{x}^{k}\right)$
then}
\item $\qquad\qquad\qquad$\textbf{set $\boldsymbol{x}^{k+1}=\mathbf{P}_{T}\boldsymbol{y}$}
\item $\qquad\qquad\qquad$\textbf{set }\emph{continue = false}
\item $\qquad\qquad$\textbf{set $\ell=\ell+1$}
\item $\qquad$\textbf{set $k=k+1$}
\end{enumerate}
Sometimes it is useful to emphasize the function $\phi$ for which
we are superiorizing, in which case we refer to the algorithm above
as the\emph{ $\phi$-superiorized version of algorithm }\textbf{\textit{$\mathbf{P}$}}\textit{.}\textit{\emph{
It is important to bear in mind that the sequence $S$ produced by
the algorithm depends also on the initial point $\bar{\boldsymbol{x}}$,
the selection of the subgradient in Line (v) of the algorithm, the
summable sequence }}$\left\{ \gamma_{\ell}\right\} _{\ell=0}^{\infty}$,
and the problem $T$. In addition, the output \emph{$O\left(T,\varepsilon,S\right)$}
of the algorithm depends on the stopping criterion $\varepsilon$.\textbf{\textit{\emph{}}}\\
\textbf{\textit{\emph{}}}\\
\textbf{\textit{\emph{Theorem 3.}}} \emph{Under the Assumptions
listed above, the Superiorized Version of Algorithm }\textbf{\textit{\emph{$\mathbf{P}$}}}\emph{
will produce a sequence $S\left(T,\bar{\boldsymbol{x}},\phi\right)$
of points in $\mathbb{R}^{J}$ that either contains a solution of
$T$ or is infinite. In the latter case, if the sequence }$\left\{ \left(\mathbf{P}_{T}\right)^{k}\boldsymbol{x}\right\} _{k=0}^{\infty}$\emph{
converges to a solution of $T$ for all $\boldsymbol{x}\in\mathbb{R}^{J}$,
then, for any $\varepsilon\in\mathbb{R}_{++}$, $O\left(T,\varepsilon,S\left(T,\bar{\boldsymbol{x}},\phi\right)\right)$
is defined and $\phi\left(O\left(T,\varepsilon,S\left(T,\bar{\boldsymbol{x}},\phi\right)\right)\right)\leq\varepsilon$.}\textbf{}\\
\textbf{Proof. }Assume that the sequence $S\left(T,\bar{\boldsymbol{x}},\phi\right)$
produced by the Superiorized Version of Algorithm \textbf{\textit{\emph{$\mathbf{P}$}}}
dos not contain a solution of $T$. We first show that in this case
the algorithm generates an infinite sequence $\left\{ \boldsymbol{x}^{k}\right\} _{k=0}^{\infty}$.
This is equivalent to saying that, for any $\boldsymbol{x}^{k}$ that
has been generated already, the condition in Line (xiii) of the algorithm
will be satisfied sooner or later (and hence $\boldsymbol{x}^{k+1}$
will be generated). This needs to happen, because as long as the condition
is not satisfied we keep resetting (in Line (xi)) the value of $\beta_{k}$
to $\gamma_{\ell}$, with ever increasing values of $\ell$. However,
$\left\{ \gamma_{\ell}\right\} _{\ell=0}^{\infty}$ is a summable
sequence of positive real numbers, and so $\gamma_{\ell}$ is guaranteed
to be arbitrarily small if $\ell$ is sufficiently large. Since $\boldsymbol{v}^{k}$
is either a unit vector in the direction of the negative subgradient
of the convex function $\phi$ at $\boldsymbol{x}^{k}$ or is the
zero vector (see Lines (v)--(viii)), \textbf{$\phi\left(\boldsymbol{x}^{k}+\beta_{k}\boldsymbol{v}^{k}\right)\leq\phi\left(\boldsymbol{x}^{k}\right)$}
must be satisfied if the positive number $\beta_{k}$ is small enough.
Also, since $\mathcal{P}r_{T}\left(\mathbf{P}_{T}\boldsymbol{x}^{k}\right)<\mathcal{P}r_{T}\left(\boldsymbol{x}^{k}\right)$
and $\mathbf{P}_{T}$ and $\mathcal{P}r_{T}$ are continuous (Assumptions
(ii) and (i), respectively), we also have that $\mathcal{P}r_{T}\left(\mathbf{P}_{T}\left(\boldsymbol{x}^{k}+\beta_{k}\boldsymbol{v}^{k}\right)\right)<\mathcal{P}r_{T}\left(\boldsymbol{x}^{k}\right)$
if $\beta_{k}$ is small enough. This completes the proof that the
condition in Line (xiii) of the algorithm will be satisfied and so
the algorithm will generate an infinite sequence $S\left(T,\bar{\boldsymbol{x}},\phi\right)$.
Observing that we have already demonstrated that the $\beta_{k}\boldsymbol{v}^{k}$
are bounded perturbations, and comparing (\ref{eq:resilient}) with
Lines (xii) and (xiv), we see that (by the bounded perturbation resilience
of $\mathbf{P}$) the assumption that the sequence $\left\{ \left(\mathbf{P}_{T}\right)^{k}\boldsymbol{x}\right\} _{k=0}^{\infty}$
converges to a solution of $T$ for all $\boldsymbol{x}\in\mathbb{R}^{J}$
implies that $S\left(T,\bar{\boldsymbol{x}},\phi\right))$ also converges
to a solution of $T$. Thus, applying Lemma 1 we obtain the final
claim of the theorem. $\square$\\

Unfortunately, this theorem does not go far enough. To demonstrate
that a methodology leads to superiorization we should be proving (under
some assumptions) a result like $\phi\left(O\left(T,\varepsilon,S\left(T,\bar{\boldsymbol{x}},\phi\right)\right)\right)<\phi\left(O\left(T,\varepsilon,R\left(T,\bar{\boldsymbol{x}}\right)\right)\right)$
in place of the weaker result at the end of the statement of the theorem.
Currently we do not have any such proofs and so we are restricted
to providing practical demonstrations that our methodology leads to
superiorization in the desired sense. In the next section we provide
such demonstrations for the Superiorized Version of Algorithm \textbf{\textit{\emph{$\mathbf{P}$}}},
for two different \textbf{\textit{\emph{$\mathbf{P}$}}}\textit{\emph{s.}}

\section{Illustrations of the superiorization methodology\label{sec:tv}}

\begin{figure}[H]
\begin{centering}
\subfloat[]{\includegraphics[scale=0.37]{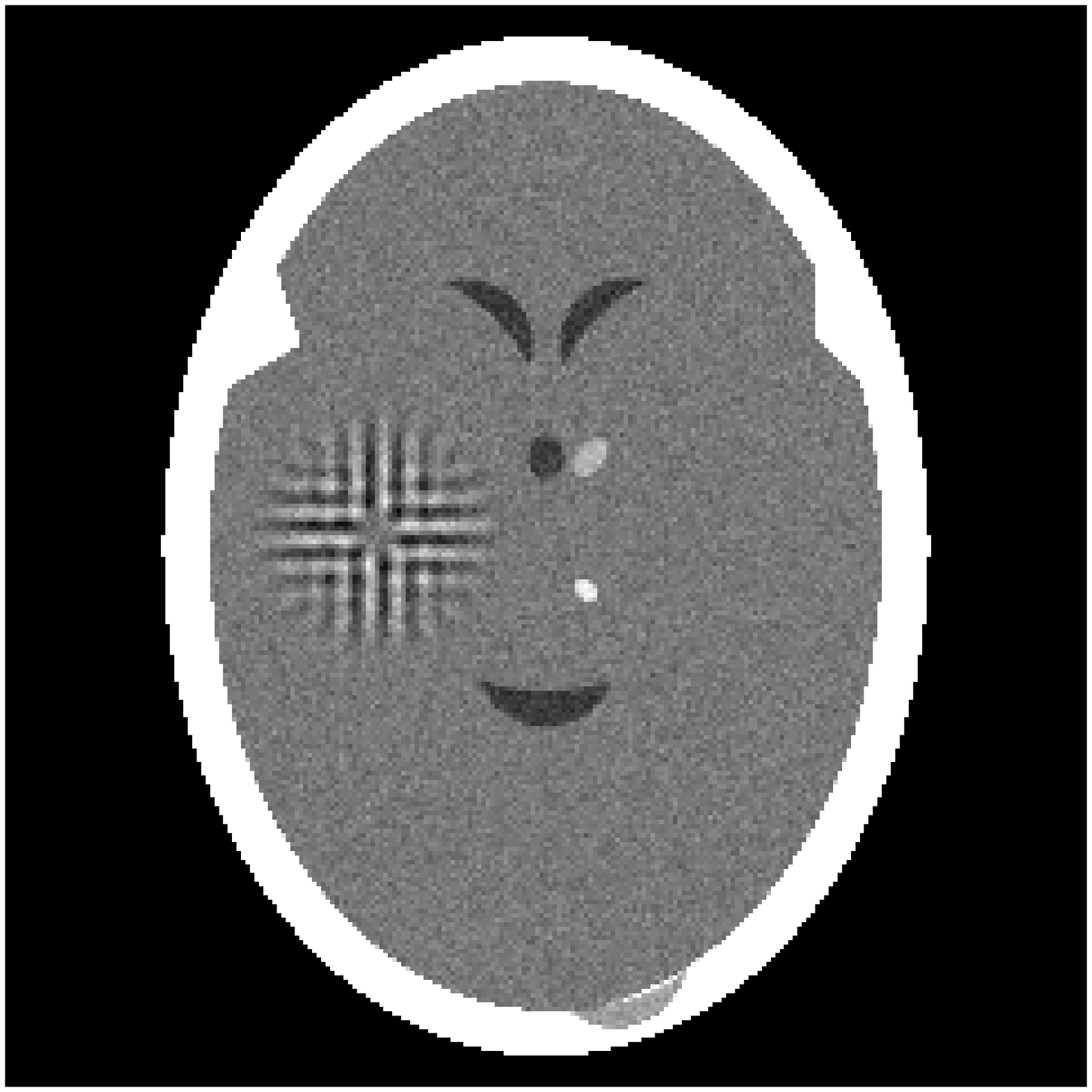}}
\par\end{centering}

\begin{centering}
\subfloat[]{

\includegraphics[scale=0.37]{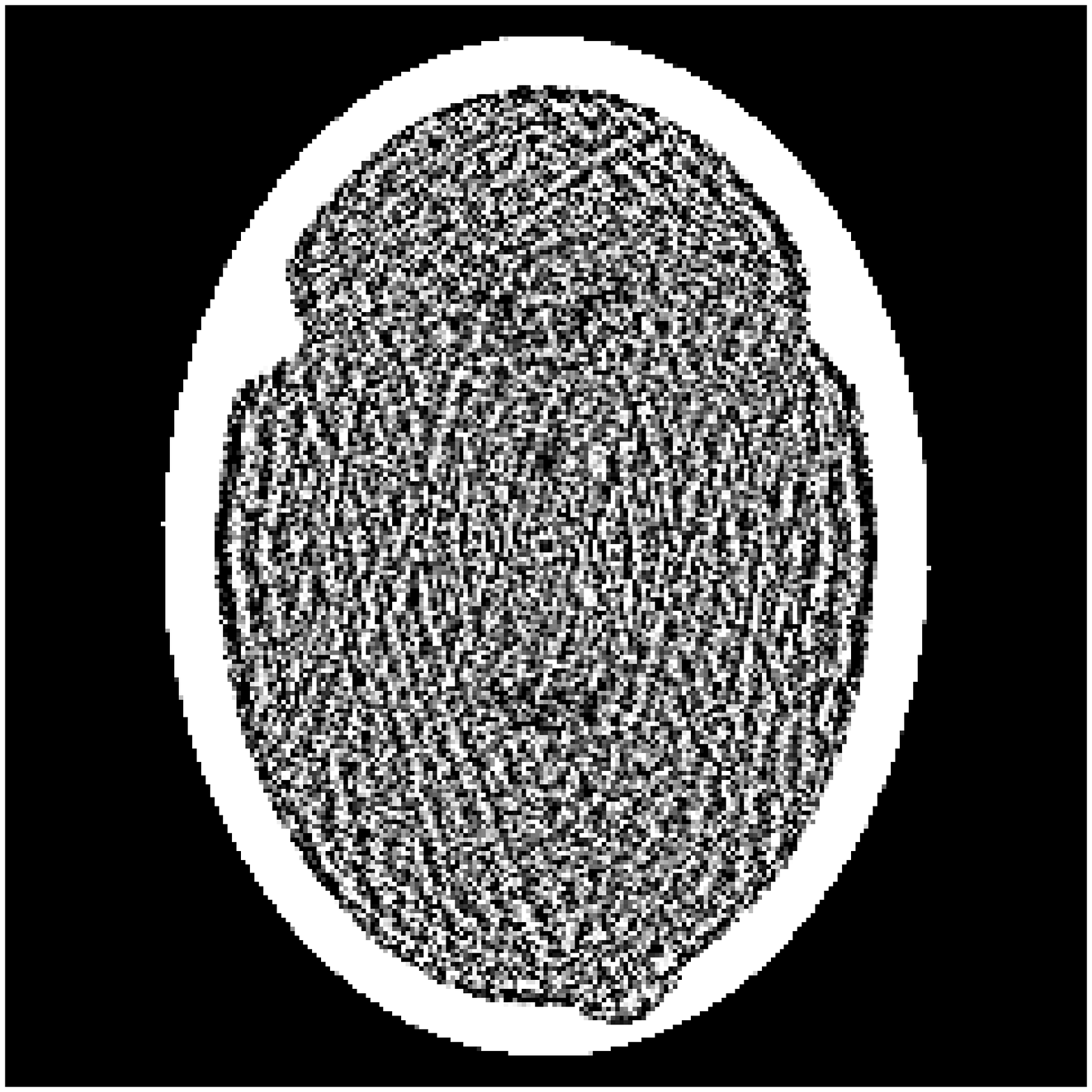}}\subfloat[]{

\includegraphics[scale=0.37]{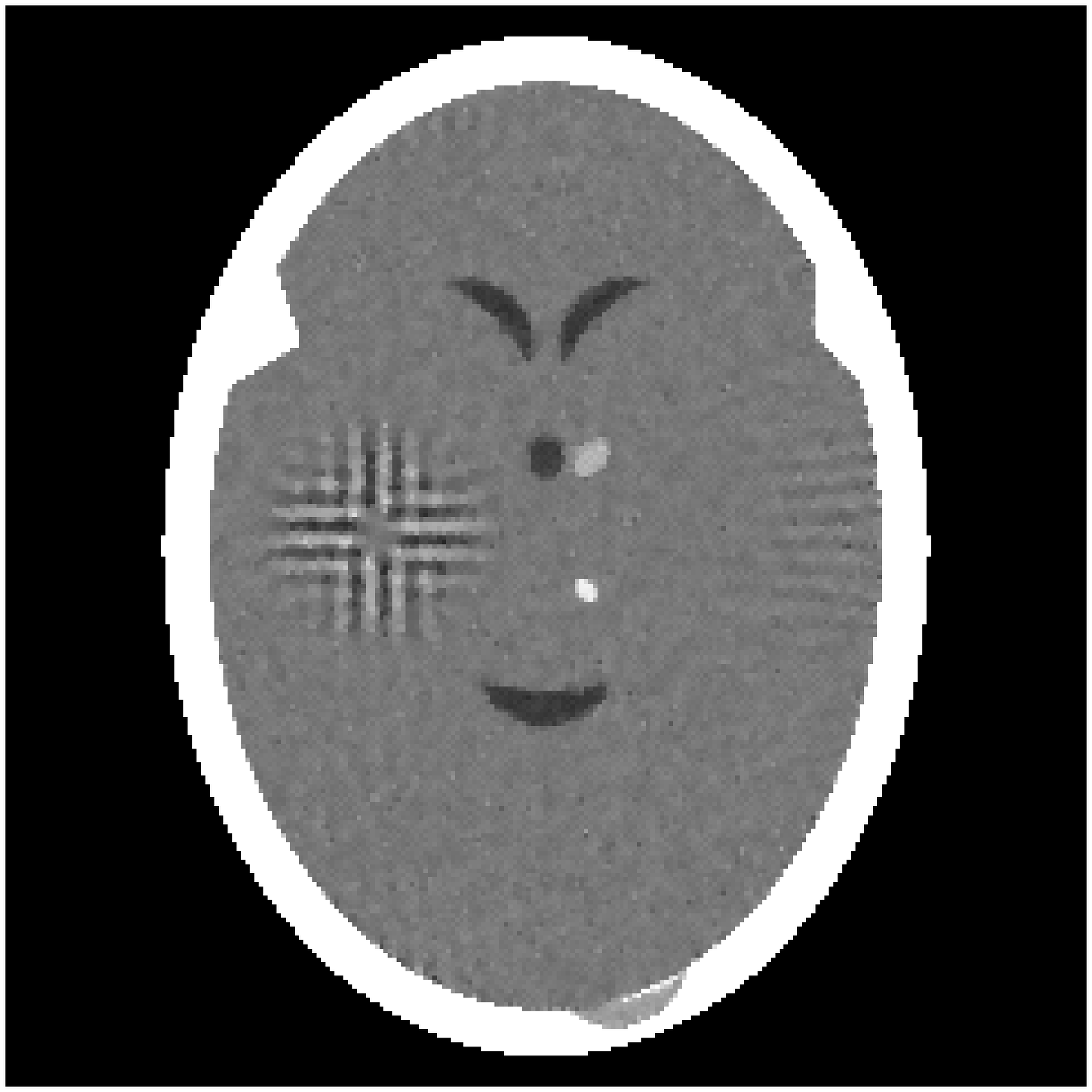}}
\par\end{centering}

\begin{centering}
\subfloat[]{

\includegraphics[scale=0.37]{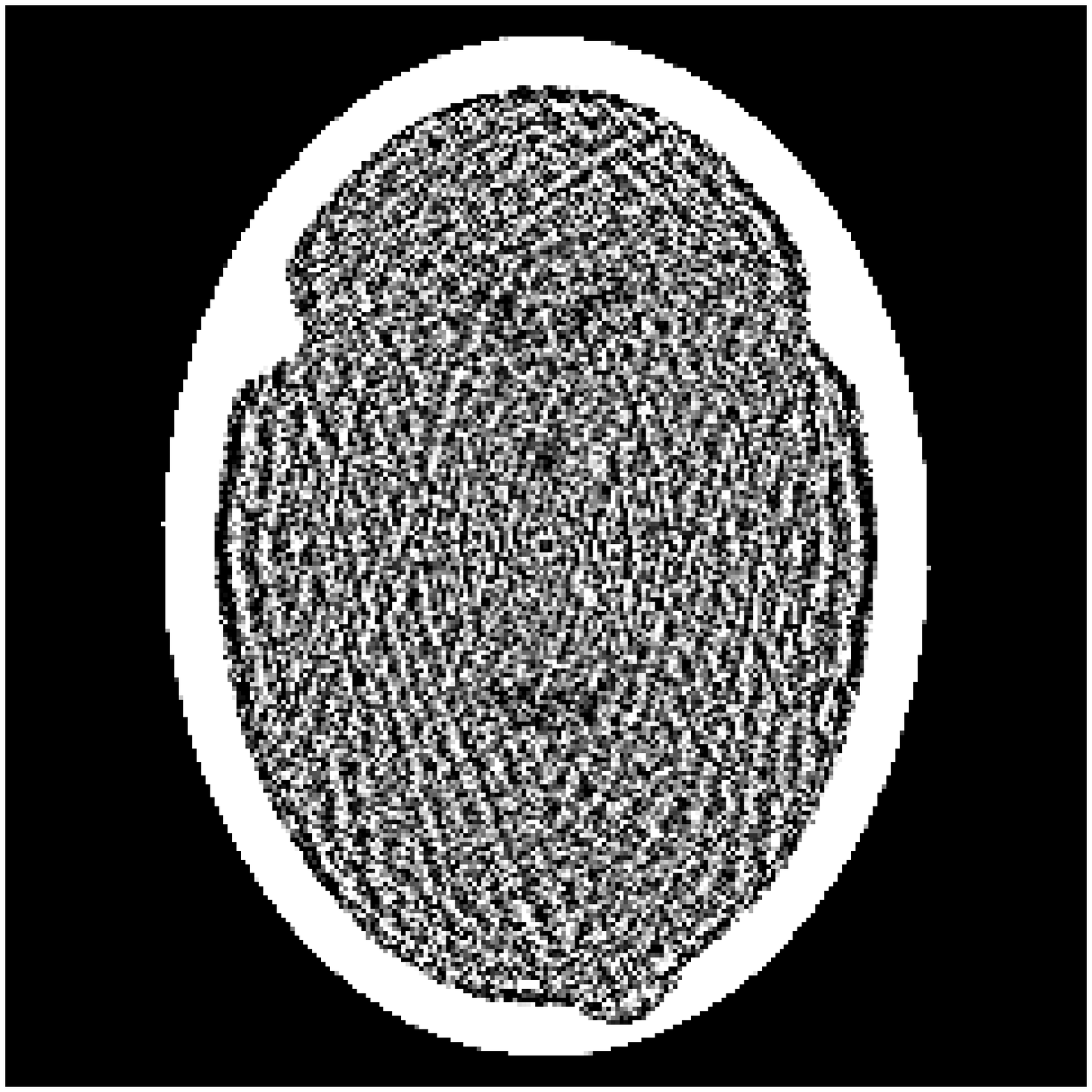}}\subfloat[]{

\includegraphics[scale=0.37]{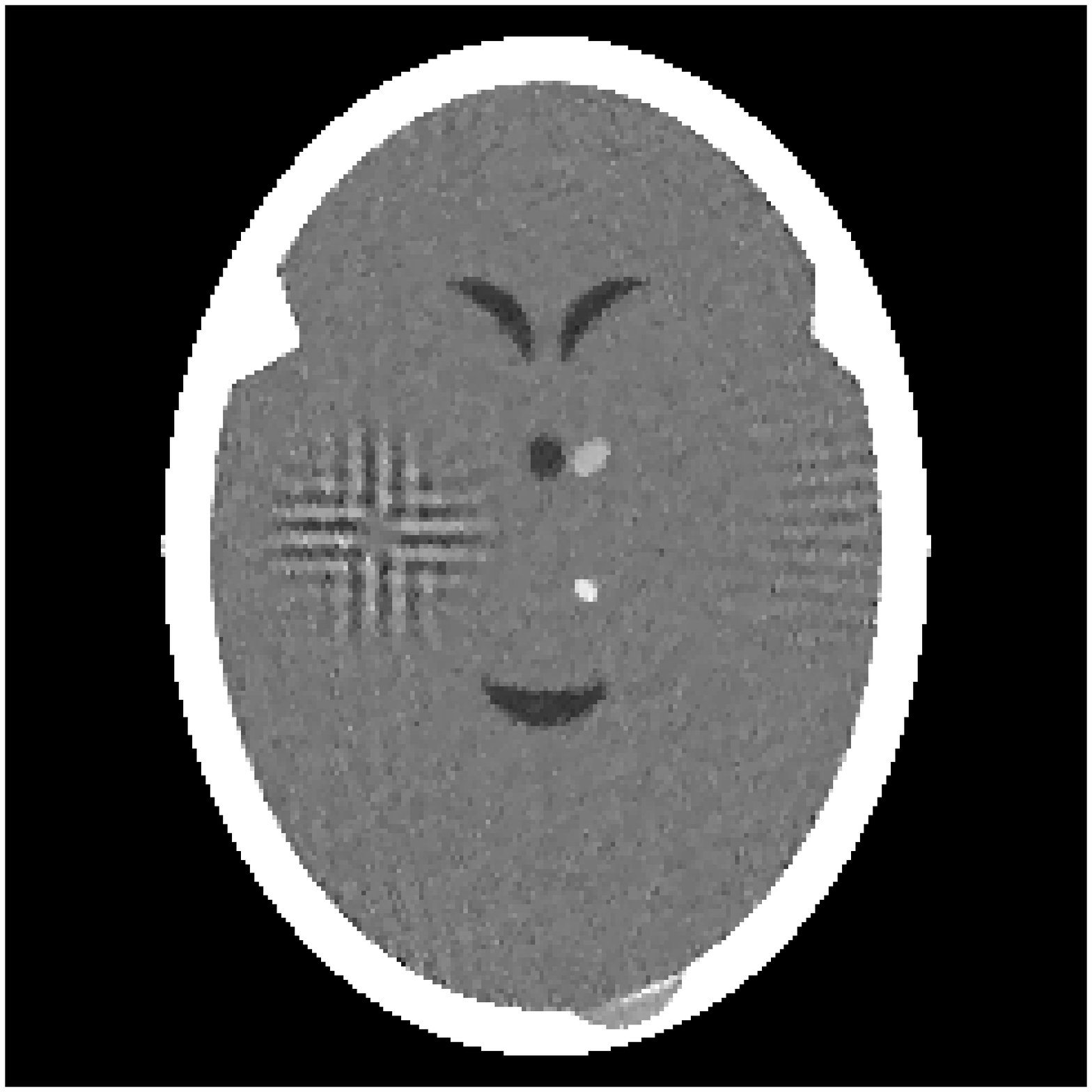}}
\par\end{centering}

\caption{\label{cap:82-consistent} A head phantom (a) and its reconstructions
from underdetermined consistent data obtained for 82~\textcolor{black}{views
using: }(b) a variant of ART, \textcolor{black}{(c) TV-superiorized
version of the same variant of ART, (d) a block-iterative projection
method, and (e) TV-superiorized version of the same block-iterative
projection method. The same initial point and stopping criterion were
used in all cases; see the text for details.}}

\end{figure}

We illustrate the superiorization methodology on a problem of reconstructing
a head cross-section (based on Figure 4.6(a) of \cite{H09}) from
its projections using both an SAP and a BIP algorithm. \textcolor{black}{(All
the computational work reported in this section was done using SNARK09
\cite{SNARK09}; the phantom, the data, the reconstructions and displays
were all generated within this same framework.)} Figure \ref{cap:82-consistent}(a)
shows a $243\times243$ digitization of the head phantom with $J=59,049$
pixels. An $\boldsymbol{x}\in\mathbb{R}^{J}$ is interpreted as a
vector of pixel values, whose components represent the average X-ray
linear attenuation coefficients (measured per centimeter) within the
$59,049$ pixels. Each pixel is of size $0.0752\times0.0752$ (measured
in centimeters). The pixel values range from $0$ to $0.5639$. For
display purposes, any value below $0.204$ is shown as black (gray
value $0$) and any value above $0.21675$ is shown as white (gray
value $255$), with a linear mapping of the pixel values into gray
values in between (the same convention is used in displaying reconstructed
images in Figures \ref{cap:82-consistent}(b)-(e)). 

\begin{table}[t]
\caption{\label{tab:Values-of-TV}Values of TV for the outputs of the various
algorithms. The second column is for the superiorized versions and
the third column is for the original versions. }
\centering\begin{tabular}{ccc}
\toprule 
Algorithm & $\phi\left(O\left(T,\varepsilon,S\left(T,\bar{\boldsymbol{x}},\phi\right)\right)\right)$ & $\phi\left(O\left(T,\varepsilon,R\left(T,\bar{\boldsymbol{x}}\right))\right)\right)$\tabularnewline
\midrule
Variant of ART & $441.50$ & $1,296.44$\tabularnewline
Variant of BIP & $444.15$ & $1,286.44$\tabularnewline
\bottomrule
\end{tabular}
\end{table}

Data were collected by calculating line integrals across the digitized
image for $82$ sets of equally spaced parallel lines, with $I=25,452$
lines in total. Each data item determines a hyperplane in $\mathbb{R}^{J}$.
Since the digitized phantom lies in the intersection of all the hyperplanes,
we have here an instance of the convex feasibility problem with a
nonempty $C$, satisfying the first condition of the statements of
Theorems 1 and 2.

For our illustration, we chose the SAP algorithm $\mathbf{P}_{\left\{ C_{1},\ldots,C_{I}\right\} }$
as determined by (\ref{eq:P_C[t]})-(\ref{eq: Bold p}) with $\Omega=\left\{ \left(1,\ldots,I\right)\right\} $
and $\omega\left(1,\ldots,I\right)=1$. This is a classical method
that in tomography would be considered a variant of the algebraic
reconstruction techniques (ART) \textcolor{black}{\cite[Chapter 11]{H09}.
For the BIP algorithm we chose $\mathbf{\mathbf{Q}_{\left\{ C_{1},\ldots,C_{I}\right\} }}$
as determined by }(\ref{eq: Blod Q operator})-(\ref{eq:definition of R})
with \textcolor{black}{$U=82$ and each block corresponding to one
of the $82$ sets of parallel lines along which the data are collected. }

The function $\phi$ for which we superiorized is defined so that,
for any $\boldsymbol{x}\in\mathbb{R}^{J}$,~$\phi\left(\boldsymbol{x}\right)$
is the \emph{total variation} (TV) of the corresponding $243\times243$
image. If the pixel values of this image are $q_{g,h}$, then the
value of the TV is defined to be \begin{equation}
\sum\limits _{g=1}^{242}\sum\limits _{h=1}^{242}\sqrt{\left(q_{g+1,h}-q_{g,h}\right)^{2}+\left(q_{g,h+1}-q_{g,h}\right)^{2}}.\label{eq:TV}\end{equation}
For the TV-Superiorized Versions of the Algorithms $\mathbf{P}_{\left\{ C_{1},\ldots,C_{I}\right\} }$
and \textcolor{black}{$\mathbf{\mathbf{Q}_{\left\{ C_{1},\ldots,C_{I}\right\} }}$}
of the previous paragraph we selected $\bar{\boldsymbol{x}}$ to be
the origin (the vector of all zeros) and $\gamma_{\ell}=0.999^{\ell}$.
Also , we \textcolor{black}{set $\varepsilon=0.01$ for the stopping
criterion, which is small compared to the }$\mathcal{P}r_{T}$\textcolor{black}{{}
of the initial point ($\mathcal{P}r_{T}\left(\bar{\boldsymbol{x}}\right)=330.208)$.}

For each of the four algorithms ($\mathbf{P}_{\left\{ C_{1},\ldots,C_{I}\right\} }$\textcolor{black}{,
$\mathbf{\mathbf{Q}}_{\left\{ C_{1},\ldots,C_{I}\right\} }$} and
their TV-superiorized versions), the sequence $S$ that is produced
by it is such that the output \emph{$O\left(T,\varepsilon,S\right)$
}is defined; see Figures \ref{cap:82-consistent}(b)-(e) for the images
that correspond to these outputs. Clearly, the superiorized reconstructions
in Figures \ref{cap:82-consistent}(c) and (e) \textcolor{black}{are
visually superior to their not superiorized versions }in Figure \ref{cap:82-consistent}\textcolor{black}{(b)
and (d), respectively. More importantly from the point of view of
our theory, consider Table \ref{tab:Values-of-TV}. As stated in the
last paragraph of the previous section, we would like to have that
}$\phi\left(O\left(T,\varepsilon,S\left(T,\bar{\boldsymbol{x}},\phi\right)\right)\right)<\phi\left(O\left(T,\varepsilon,R\left(T,\bar{\boldsymbol{x}}\right)\right)\right)$.
While we are not able to prove that this is the case in general, \textcolor{black}{Table
\ref{tab:Values-of-TV} clearly shows it to be the case for the two
algorithms discussed in this section.}

A final important point that is illustrated by the experiments in
this section is that, from the practical point of view, TV-superiorization
is as useful as TV-optimization. This is because a realistic phantom,
such as the one in Figure \ref{cap:82-consistent}(a), is unlikely
to be TV-minimizing subject to the constraints provided by the measurements.
In fact, the TV value of our phantom is $450.53$, which is larger
than that for either of the TV-superiorized reconstructions in the
second column of Table \textcolor{black}{\ref{tab:Values-of-TV}.
While an optimization method should be able to find an image with
a lower TV value, there is no practical point for doing that. The
underlying aim of what we are doing is to estimate the phantom from
the data and producing an image whose TV value is further from the
TV value of the phantom than that of our superiorized reconstructions
is unlikely to be helpful towards achieving this aim.}

\section{Discussion and conclusions\label{sec:conclusion}}

Stability of algorithms under perturbations is generally studied in
numerical analysis with the aim of proving that an algorithm is stable
so that it can {}``endure'' all kinds of imperfections in the data
or in the computational performance. Here we have taken a proactive
approach designed to extract specific benefits from the kind of stability
that we term perturbation resilience. We have been able to do this
in a context that includes, but is much more general than, feasibility-optimization
for intersections of convex sets. 

Our premise has been that (1) there is available a bounded perturbations
resilient iterative algorithm that solves efficiently certain type
of problems and (2) we desire to make use of perturbations to find
for these problems\emph{ }solutions that,\emph{ }according to some
criterion, are superior to the ones to which we would get without
employing perturbations. To accomplish this one must have a way of
introducing perturbations that take into account the criterion according
to which we wish to {}``superiorize'' the solutions of the problems.

We have set forth the fundamental principle, have given some mathematical
formulations and results, and have shown potential benefits (in the
field of image reconstruction from projections). However, the superiorization
methodology needs to be studied further from the mathematical, algorithmic
and computational points of view in order to unveil its general applicability
to inverse problems. As algorithms are developed and tested a dialog
on algorithmic developments must be accompanied by mathematical validation
and applications to simulated and real data from various relevant
fields of applications.

Validating the concept means proving precise statements about the
behavior of iterates $\left\{ \boldsymbol{x}^{k}\right\} _{k=0}^{\infty}$
generated by the superiorized versions of algorithms. Under what conditions
do they converge? Can their limit points be characterized? How would
different choices of the perturbation coefficients $\beta_{k}$ and
the perturbation vectors $\boldsymbol{v}^{k}$ affect the superiorization
process? Can different schemes for generating the $\beta_{k}$s be
developed, implemented, investigated? Enlarging the arsenal of bounded
perturbation resilience algorithms means generalizing existing proofs
of bounded perturbations resiliency of algorithms and developing new
theories that will bring more algorithms into the family of bounded
perturbations resilient algorithms. Further developments should include
the problem of finding a common fixed point of a family of operators
(a direct generalization of the convex feasibility problem) and studying
the behavior of superiorization algorithms in inconsistent situations
when the underlying solution set is empty. Thus we view the material
in this paper as only an initial step in a promising new field of
endeavor for solving inverse problems.

\ack{}{}

This work was supported by Award Number R01HL070472 from the National
Heart, Lung, And Blood Institute. The content is solely the responsibility
of the authors and does not necessarily represent the official views
of the National Heart, Lung, And Blood Institute or the National Institutes
of Health.


\section*{References}
\begin{thebibliography}{33}
\bibitem{ac89}Aharoni R and Censor Y 1989 Block-iterative projection
methods for parallel computation of solutions to convex feasibility
problems \emph{Linear Algebra Appl.} \textbf{120}, 165--75 

\bibitem{bb96} Bauschke HH and Borwein JM 1996 On projection algorithms
for solving convex feasibility problems \emph{SIAM Rev.} \textbf{38}
367\textendash{}426

\bibitem{Baus06} Bauschke HH, Combettes PL and Kruk SG 2006 Extrapolation
algorithm for affine-convex feasibility problems\emph{ Numer. Algorithms}
\textbf{41} 239--74

\bibitem{bmr03} Bauschke HH, Matoušková E and Reich S 2004 Projection
and proximal point methods: convergence results and counterexamples.
\emph{Nonlinear Anal.} \textbf{56} 715--38 

\bibitem{Butn90} Butnariu D and Censor Y 1990 On the behavior of
a block-iterative projection method for solving convex feasibility
problems \emph{Int. J. Comput Math.} \textbf{34} 79--94

\bibitem{bdhk07}Butnariu D, Davidi R, Herman GT and Kazantsev IG
2007 Stable convergence behavior under summable perturbations of a
class of projection methods for convex feasibility and optimization
problems \emph{IEEE J. Sel. Top. Sign. Process.} \textbf{1} 540--7

\bibitem{brz05}Butnariu D, Reich S and Zaslavski AJ 2006 Convergence
to fixed points of inexact orbits of Bregman-monotone and nonexpansive
operators in Banach spaces \emph{Fixed Point Theory and Applications
}ed H F Nathansky, B G de Buen, K Goebel, W A Kirk and B Sims (Yokohama:
Yokohama Publishers) pp 11\textendash{}32

\bibitem{brz08} Butnariu D, Reich S and Zaslavski AJ 2008 Stable
convergence theorems for infinite products and powers of nonexpansive
mappings \emph{Numer. Func. Anal. Opt.} \textbf{29} 304--23

\bibitem{byrnebook} Byrne CL 2008\emph{ Applied Iterative Methods}
(AK Peters)

\bibitem{cccdh09}Censor Y, Chen W, Combettes PL, Davidi R and Herman
GT 2009 On the effectiveness of projection methods for convex feasibility
problems with linear inequality constraints, \emph{Opt. Online}, http://www.optimization-online.org/DB\_HTML/2009/12/2500.html

\bibitem{ceh01}Censor Y, Elfving T and Herman GT 2001 Averaging strings
of sequential iterations for convex feasibility problems. \emph{Inherently
Parallel Algorithms in Feasibility and Optimization and Their Applications}
ed Butnariu D, Censor Y and Reich S (Elsevier Science Publishers)
pp 101--14

\bibitem{cgg01} Censor Y, Gordon D and Gordon R 2001 BICAV: A block-iterative,
parallel algorithm for sparse systems with pixel-related weig\textcolor{black}{hting
}\textcolor{black}{\emph{IEEE Trans. Med. Imaging}}\textcolor{black}{{}
}\textbf{\textcolor{black}{20}} 1050--60

\bibitem{cs09}Censor Y and Segal A 2009 On the string averaging method
for sparse common fixed points problems \textcolor{black}{\emph{Int.
Trans. Oper. Res.}}\textcolor{black}{{} }\textbf{16} 481--94 

\bibitem{ct03}Censor Y and Tom E 2003 Convergence of string-averaging
projection schemes for inconsistent convex feasibility problems \emph{Optim.
Methods Softw.} \textbf{18} 543--54

\bibitem{CZ97} Censor Y and Zenios \textit{\emph{SA 1997 }}\textit{Parallel
Optimization: Theory, Algorithms and Applications} (Oxford University
Press)

\bibitem{chinneck-book} Chinneck JW 2007 \emph{Feasibility and Infeasibility
in Optimization: Algorithms and Computational Methods} (Springer)

\bibitem{c96} Combettes PL 1996 The convex feasibility problem in
image recovery\textcolor{black}{{} }\textcolor{black}{\emph{Adv. Imag.
Elec. Phys.}}\textcolor{black}{{} }\textbf{\textcolor{black}{95}}\textcolor{black}{{}
}155--270 

\bibitem{Jamo97} Combettes PL 1997 Hilbertian convex feasibility
problem: Convergence of projection methods \emph{Appl. Math. Opt.}
\textbf{35} 311--30 

\bibitem{Imag97} Combettes PL 1997 Convex set theoretic image recovery
by extrapolated iterations of parallel subgradient projections \emph{IEEE
T. Image Process.} \textbf{6} 493--506

\bibitem{tv2}Combettes PL and Luo J 2002 An adaptive level set method
for nondifferentiable constrained image recovery \emph{IEEE T. Image
Process.} \textbf{11} 1295--304 

\bibitem{tv3}Combettes PL and Pesquet JC 2004 Image restoration subject
to a total variation constraint \emph{IEEE T. Image Process.}\textbf{
13} 1213--22

\bibitem{crombez02} Crombez G 2002 Finding common fixed points of
strict paracontractions by averaging strings of sequential iterations
\emph{J. Nonlinear Convex Anal} \textbf{3} 345--51

\bibitem{dhc09} Davidi R, Herman GT and Censor Y 2009 Perturbation-resilient
block-iterative projection methods with application to image reconstruction
from projection\textcolor{black}{s }\textcolor{black}{\emph{Int. Trans.
Oper. Res.}}\textcolor{black}{{} }\textbf{16} 505--24

\bibitem{SNARK09} Davidi R, Herman GT and Klukowska J 2009\emph{
}SNARK09: A programming system for the reconstruction of 2D images
from 1D projections (http://www.snark09.com/)

\bibitem{eGG81} Eggermont PPB, Herman GT and Lent A 1981 Iterative
algorithms for large partitioned linear systems, with applications
to image reconstruction \emph{Linear Algebra Appl.} \textbf{40} 37--67 

\bibitem{Gonz01} González-Castaño FJ, García-Palomares UM, Alba-Castro
JL and Pousada-Carballo JM 2001 Fast image recovery using dynamic
load balancing in parallel architectures, by means of incomplete projections
\emph{IEEE T. Image Process.} \textbf{10} 493--99 

\bibitem{H09} Herman GT 2009\emph{ }\textit{Fundamentals of Computerized
Tomography: Image Reconstruction from Projections}\emph{ }2nd ed.
(Springer)

\bibitem{hd08}Herman GT and Davidi R 2008 On image reconstruction
from a small number of projections\emph{ Inverse Problems} \textbf{24}:045011

\bibitem{Lopu97} Kiwiel KC and \L{}opuch B 1997 Surrogate projection
methods for finding fixed points of firmly nonexpansive mappings \emph{SIAM
J. Optim.} \textbf{7} 1084--1102

\bibitem{Otta88} Ottavy N 1988 Strong convergence of projection-like
methods in Hilbert spaces \emph{J. Optim. Theory Appl. }\textbf{56}
433--461

\bibitem{pen09} Penfold SN, Schulte RW, Censor Y, Bashkirov V, McAllister
S, Schubert KE, Rosenfeld AB (to appear) Block-iterative and string-averaging
projection algorithms in proton computed tomography image reconstruction
ed Censor Y, Jiang M and Wang G \emph{Biomedical Mathematics: Promising
Directions in Imaging, Therapy Planning and Inverse Problems} (Medical
Physics Publishing) 

\bibitem{Pier84}Pierra G 1984 Decomposition through formalization
in a product space \emph{Math. Program.} \textbf{28} 96--115

\bibitem{rh03}Rhee H 2003 An application of the string averaging
method to one-sided best simultaneous approximation \emph{J. Korea
Soc. Math. Educ. Ser. B Pure Appl. Math. }\textbf{10} 49--56
\end{thebibliography}
\end{document}